\documentclass{amsart}

\usepackage{amssymb}
\usepackage[all]{xy}

\newtheorem{thm}{Theorem}

\newtheorem{lem}[thm]{Lemma}
\newtheorem{cor}[thm]{Corollary}

\newtheorem{prop}[thm]{Proposition}

\theoremstyle{definition}
\newtheorem{defn}[thm]{Definition}

\newtheorem{say}[thm]{}
\newtheorem{exmp}[thm]{Example}

\newtheorem{ques}[thm]{Question}    

\newtheorem{rem}[thm]{Remark}          

\newtheorem*{ack}{Acknowledgments}      

\newtheorem{defn-thm}[thm]{Definition--Theorem}  
\newtheorem{defn-lem}[thm]{Definition--Lemma}  

\newtheorem{comments}[thm]{Comments}

\theoremstyle{remark}

\let \cedilla =\c
\renewcommand{\c}[0]{{\mathbb C}}  

\renewcommand{\o}[0]{{\mathcal O}} 

\renewcommand{\a}[0]{{\mathbb A}}

\newcommand{\qtq}[1]{\quad\mbox{#1}\quad}
\newcommand{\spec}[0]{\operatorname{Spec}}

\newcommand{\supp}[0]{\operatorname{Supp}}    
\newcommand{\red}[0]{\operatorname{red}}    
\newcommand{\codim}[0]{\operatorname{codim}}

\newcommand{\chr}[0]{\operatorname{char}}

\newcommand{\depth}[0]{\operatorname{depth}}

\newcommand{\cond}[0]{\operatorname{cond}}

\def\into{\DOTSB\lhook\joinrel\to}

\def\loccoh#1.#2.#3.#4.{H^{#1}_{#2}(#3,#4)}

\DeclareMathAlphabet{\mathchanc}{OT1}{pzc}%
                                {m}{it}

\usepackage[all]{xy}\xyoption{dvips}

\begin{document}
\bibliographystyle{amsalpha}


\title{Variants of normality for   Noetherian schemes}
\author{J\'anos Koll\'ar}

\begin{abstract} 
This note presents a  uniform treatment of normality and three of its 
variants---topological, weak and seminormality---for
Noetherian schemes. The key is to define these notions for pairs
$(Z, X)$ consisting  of a (not necessarily reduced) scheme $X$ and a  closed, 
nowhere dense subscheme $Z$. 
An advantage of the new definitions is that, unlike the usual absolute ones, 
 they are preserved by completions.
This shortens some of the proofs and leads to more general  results.
\end{abstract}

\maketitle

\begin{defn}  \label{fin.mod.defn}
Let $X$ be a  scheme and $Z\subset  X$ a closed, nowhere dense subscheme.
A {\it finite modification} of $X$ centered at $Z$ is a 
finite morphism $p:Y\to X$ such that none of the associated primes of $Y$ is contained in $Z_Y:=p^{-1}(Z)$ and 
$$
p|_{Y\setminus Z_Y}:Y\setminus Z_Y\to X\setminus Z
 \qtq{is an isomorphism.}
\eqno{(\ref{fin.mod.defn}.1)}
$$
Let $j:X\setminus Z\into X$ be the natural injection and $J_Z\subset \o_X$
the largest subsheaf supported on $Z$. There is a one-to-one
correspondence between finite modifications and coherent $\o_X$-algebras
$$
\o_X/J_Z\subset p_*\o_Y\subset j_*\o_{X\setminus Z}.
\eqno{(\ref{fin.mod.defn}.2)}
$$
The notion of an {\it integral modification} of $X$ centered at $Z$
is defined analogously.
Let $A\subset j_*\o_{X\setminus Z}$ be the largest subalgebra that is
integral over $\o_X$. Then $\spec_XA$ is the maximal  integral modification, called the 
{\it relative normalization} of  the pair $Z\subset X$. We denote it  by
$$
\pi: (Z^{\rm rn}\subset X^{\rm rn})\to (Z\subset X)\qtq{or by}
\pi: X_Z^{\rm rn}\to X.
\eqno{(\ref{fin.mod.defn}.3)}
$$
The relative normalization is the limit of all  finite modifications
centered at $Z$.
(It would be called the relative normalization of $X$ in $X\setminus Z$  in the terminology of \cite[Tag 0BAK]{stacks-project}.)

If $(x, X)$ is semilocal  then a finite modification of $X$ centered at
 $Z=\{x\}$ is called a 
 {\it punctual modification} of $(x,X)$.
The relative normalization of $(x, X)$ is called the
{\it punctual normalization} and denoted by
$\pi: \bigl(x^{\rm pn}, X^{\rm pn}\bigr)\to (x, X)$.
\end{defn}

\begin{defn}\label{5.props.den}
 Let $X$ be a Noetherian scheme and $Z\subset  X$ a closed, nowhere dense subscheme. We define 4 properties of such pairs, depending on the 
behavior of all finite modifications  $p:Y\to X$  centered at $Z$.
  
\begin{enumerate}
\item[N] (normality): every  $p:Y\to X$ is an isomorphism.
\item[TN]  (topological normality):   every  $p:Y\to X$ is a universal homeomorphism (Definition \ref{univ.homeo.defn}).
\item[WN]  (weak normality):   if $p:Y\to X$ is a universal homeomorphism then
it is an isomorphism.
\item[SN] (seminormality):   if $p:Y\to X$ is a universal homeomorphism that preserves residue fields (Definition \ref{univ.homeo.defn}) then
it is an isomorphism.
\end{enumerate}
We stress that  $X$ is an arbitrary Noetherian scheme and 
these are all properties of pairs $Z\subset  X$. 
A direct predecessor of these definitions is in the works of S.~Kov\'acs who developed the notions of
rational  and Du~Bois pairs \cite{MR2784747, Kovacs11a}, but
 working with pairs instead of  schemes has long been a
 theme of
the Minimal Model Program; see  \cite{k-pairs, laz-book, k-db, kk-singbook}.

We also say that 
``$Z\subset X$ is  a normal pair''  and similarly for the  other variants.
A semilocal scheme $(x, X)$ is called {\it punctually normal}
if the pair $\{x\}\subset X$ is normal; similarly for the other versions.

We say that a scheme $X$ without isolated points satisfies N (or TN, ...) 
if the pair $Z\subset X$ satisfies N (or TN, ...) 
for every closed, nowhere dense subscheme $Z$. If $X$ has isolated points
then in addition we require these points to be reduced for N, WN and SN. 

It would have been possible to formulate all  these notions for
integral morphisms, but I prefer to stay with Noetherian schemes.
In all results the integral morphism versions are direct consequences of the
finite morphism versions.
\end{defn}

In the definition of topological normality it seems artificial to
restrict to maps that are isomorphisms over $X\setminus Z$; one should
clearly allow universal homeomorphisms over $X\setminus Z$. This leads to the following variant of TN.

\begin{defn}  \label{5.props.den.1}
Let $X$ be a  scheme and $Z\subset  X$ a closed, nowhere dense subscheme.
A {\it finite topological modification} of $X$ centered at $Z$ is a 
finite morphism $p:Y\to X$ such that none of the associated primes of $Y$ is contained in $Z_Y:=p^{-1}(Z)$ and 
$$
p|_{Y\setminus Z_Y}:Y\setminus Z_Y\to X\setminus Z
 \qtq{is a universal homeomorphism.}
\eqno{(\ref{5.props.den.1}.1)}
$$
In analogy with the notion of topological normality we introduce the
following.
\begin{enumerate}
\item[STN]  (strong topological normality): every finite topological modification of $X$ centered at $Z$ is a universal homeomorphism.
\end{enumerate}
This is clearly a  topological property. That is, let
$f:Y\to X$ be a  finite, universal homeomorphism. Then
$Z\subset X$  is STN iff
$f^{-1}(Z)\subset Y$ is STN.
In particular, $X$ is STN iff $\red X$ is STN.

It is clear that STN $\Rightarrow$ TN; 
eventually we show in Theorem \ref{TN=STN.factor.cor} that STN  is equivalent to TN.

\end{defn}

\begin{comments}\label{basic.def.comms} The main contribution of this note
is Definitions \ref{5.props.den}--\ref{5.props.den.1}, or rather the assertion the 
these concepts are natural and useful.  To support this claim, we show the following. 
\begin{itemize}
\item  A scheme satisfies N (resp.\ WN or SN) 
iff it is normal (resp.\ weakly or seminormal) (Proposition \ref{pf.of.hrt.defn.2} and Definition \ref{weak.norm.defn}).
\item  All 5  properties in Definitions \ref{5.props.den}--\ref{5.props.den.1} are local, even eqivalent to the punctual versions (Proposition \ref{local.nature.prop}).
\item All the  properties are preserved by completion (Corollary \ref{formal.nature.thm}).
\item A scheme $X$ is normal (resp.\ weakly or seminormal) iff the completion  $(\hat x, \hat X)$ is punctually normal (resp.\ weakly or seminormal) for every $x\in X$ (Corollary \ref{n.w.s.punctual.cor}). 
\item All the properties descend for  faithfully flat  morphisms (Proposition \ref{hart.ff.up.lem}).
\item All the properties ascend for  regular morphisms (Theorem \ref{reduced.morph.preserves.P}).
\item Properties TN,  STN and being geometrically unibranch coincide (Theorem \ref{TN=STN.factor.cor}
 and Corollary \ref{TN=unibranch.cor}).
\item Weak and seminormality can be described in terms of the conductors of the punctual normalizations (Corollary \ref{sn.wn.loc.char.cor}).  
\item Finiteness of relative normalization is a punctual property (Theorem \ref{punc.nag.lem}).
\item  Nagata schemes can be  characterized using punctual normalizations or
the reducedness of the formal fibers (Theorem \ref{nagata.char.thm}).
\end{itemize}
Note, however, that many of the  arguments are either classical or can be traced back to one of the main sources
\cite{ EGA, MR0142547, sga2, MR0277542, MR571055, MR556308}; I try to give more precise references in each section. 
For many of the topics the best reference is \cite{stacks-project}.

(\ref{basic.def.comms}.1) It would be natural to say that  ``$(x, X)$ is normal'' or ``$X$  is normal along $Z$''
but this is at variance with standard usage since $Z\subset  X$ can be  a normal pair even if $X$ is not normal at any point of $Z$. 
For instance, the pairs
$$
\{(0,0,0)\}\subset X_1:=(xy=0)\subset \a^3\qtq{and} \{(0,0,0)\}\subset X_2:=(z^2=0)\subset \a^3
$$
are both normal by our definition. This may seem somewhat perverse, but the resulting flexibility is a crucial ingredient in several of our proofs.

(\ref{basic.def.comms}.2)  A key advantage of the above definitions is that the punctual versions are preserved by completion (\ref{formal.nature.thm}). 
This makes it possible to study several basic results that were known for
 Nagata schemes and extend them to general Noetherian schemes.  In addition, the resulting proofs are  shorter and use  simpler inductive steps.
Note that completions of normal Noetherian rings can be quite complicated;
see the examples in \cite[A.1]{MR0155856} and 
 \cite{MR846452,MR1102888,MR2068083,MR2892767}.

(\ref{basic.def.comms}.3) It would seem to be natural to introduce another notion  TNR   (topological normality preserving residue fields):   every finite modification $p:Y\to X$ centered at $Z$ is a universal homeomorphism and preserves residue fields.
However,   example (\ref{tnr.not.top.exmp}) shows that TNR is not topological. Thus, although
the introduction of TNR makes the diagram in (\ref{basic.def.comms}.4)
pleasingly symmetrical, it is probably not a useful concept.

(\ref{basic.def.comms}.4) For any pair $Z\subset  X$ we have the following obvious implications
between these notions
$$
\begin{array}{ccccc}
TN & + & WN & 
\Leftrightarrow & N\\
\Uparrow &&\Downarrow  &&\\
TNR & + & SN & \Leftrightarrow& N.
\end{array}
$$
\end{comments}

\begin{ack} I thank A.~J.~de~Jong, W.~Heinzer, R.~Heitmann, S.~Kov\'acs
and S.~Loepp   for helpful
comments and references.
Partial financial support    was provided  by  the NSF under grant number
 DMS-1362960.
\end{ack}

\section{Normal pairs}

First we show that our definition of normality coincides with the usual one.

\begin{prop}\label{pf.of.hrt.defn.2} 
A Noetherian scheme $X$ without isolated points 
is normal (resp.\ weakly or seminormal) iff  $Z\subset X$ is a normal  (resp.\ weakly or seminormal) pair
for every nowhere dense closed subscheme $Z$
\end{prop}

Proof. Assume that $X$ is normal. Every finite modification centered
at $Z$ is dominated by the normalization, hence necessarily an isomorphism.
Similarly, if $X$ is weakly (resp.\ semi) normal then every birational, universal homeomorphism (resp.\ one that also preserves residue fields)
is an isomorphism.

To see the converse, we may assume that $X$ is affine (cf.\ Proposition \ref{local.nature.prop}).
First we show that $X$ is reduced. If not then 
there is an  irreducible subscheme $Z$ such that $J$, the ideal of all
nilpotent sections with support in $Z$, is nonzero.
If there is such a $Z$ that is nowhere dense 
then   $\spec_X\o_X/J\to X$ shows that  $Z\subset X$
is not a seminormal pair. Otherwise $X$ has no embedded
points hence there is an $r\in \o_X$ such that
multiplication by $r$ is injective but not surjective on $J$. 
(Here we use that $\dim Z>0$.)
There is a natural algebra structure on the $\o_X$-module
$\o_X+\tfrac1{r}J$ where we set $\tfrac1{r}J\cdot \tfrac1{r}J=0$.
The diagonally embedded   $\delta:J\into \o_X+\tfrac1{r}J$
is an ideal and  $\spec_X\bigl((\o_X+\tfrac1{r}J\bigr)/\delta(J)\bigr)\to X$
shows that  $\bigl((r=0)\cap Z\subset X\bigr)$ is not a seminormal pair.
 Thus $X$ is reduced.

Let $\bar X\to X$ be the normalization (resp.\ weak or seminormalization).
We are done if $X=\bar X$. Otherwise pick any $\phi\in \o_{\bar X}\setminus \o_X$
and set $Z:=\supp \bigl(\o_X[\phi]/\o_X\bigr)$.
Then $\spec_X\o_X[\phi]\to X$ shows that  the pair 
$Z\subset X$ is not  normal (resp.\ weakly or seminormal).\qed

\medskip

\begin{lem}\label{hart.serre.loc.thm} Let $(R,m)$ be a local Noetherian ring. Then exactly one of the following holds.
\begin{enumerate}
\item $(R,m)$ is Artinian, 
\item  $(R,m)$ is regular of dimension $ 1$, 
\item  $\depth_mR\geq 2$ or
\item $(R,m)$ is a non-normal pair.
\end{enumerate}
\end{lem}

Note that $(R,m)$ is a non-normal pair iff  there is a ring homomorphism
$\phi:R\to S$ (other than an isomorphism) whose kernel and cokernel are killed by a power of $m$ and such that $m$ is not an associated prime of $S$.

\medskip

Proof.  The  proof is essentially the same as that of Serre's criterion of
 normality.  
 $(R,m)$ is not Artinian iff $V(m)\subset \spec R$ is nowhere dense; we assume this from now on.

Let $J\subset R$ be the largest ideal killed by a power of $m$.
If $J\neq 0$ then   $R\to R/J$ shows that $(R,m)$ is  a non-normal pair.

Otherwise $J=0$ and there is an  $r\in m$ that is not a zero-divisor. If
$m$ is not an associated prime of $R/(r)$ then  $\depth_mR\geq 2$.
Thus we are left with the case when there is an $a\in R\setminus (r)$ such that
$am\subset rR$.  

If $am\subset rm$ then, by the determinantal trick \cite[Prop.2.4]{at-mc}, 
  $a/r$ satisfies a monic polynomial,
hence  $R\subset R[a/r]$ shows that   $(R,m)$ is  a non-normal pair.

Otherwise  there is a $t_0\in m$ such that  $at_0=r$; in particular
$a$ is not a zero-divisor.
For any $t\in m$ we have  $at=rt'$ for some $t'\in R$. 
Thus $a(t-t't_0)=at-t'(at_0)=at-t'r=0$.  
Since $a$ is not a zero-divisor this implies that 
 $t=t't_0$ and so
$m=(t_0)$.  Thus  $(R,m)$ is regular of dimension 1. \qed
\medskip

This directly  implies the following  claims.


\begin{cor} \label{cart.pair.nomal.cor}
Let $X$ be a Noetherian scheme and $D\subset X$ an effective 
Cartier divisor. The pair $D\subset X$ is  normal iff $X$ is regular at the 
generic points of $D$ and  $D$ has no embedded points. \qed
\end{cor}

\begin{cor} A Noetherian scheme $X$ is $S_2$ iff $Z\subset X$ is a normal pair
for every closed subscheme of codimension $\geq 2$. \qed
\end{cor}

\begin{say}\label{props.from.relmorm.say.1} 
All of the  properties in Definition \ref{5.props.den} can be understood in terms of the
relative normalization $\pi_Z: X_Z^{\rm rn}\to X$. 

For TN  this is clear.
 Since $\pi_Z: X_Z^{\rm rn}\to X$ is the limit of finite modifications, 
$Z\subset X$ satisfies TN   iff $\pi$  is
a  universal homeomorphism.
 The only complication is that the
relative normalization need not be Noetherian in general. Luckily
we will be able to avoid this issue.

A description of  WN and SN in terms of the conductor of the relative normalization is given in (\ref{sn.wn.loc.char.cor}); this is less straightforward.
\end{say}

\section{Punctual nature}

For a scheme $X$ and point $x\in X$ we set $X_x:=\spec_X\o_{x,X}$. 
Thus $(x, X_x)$ is a local scheme. 
We show that our properties  can be described in terms of the local schemes of $X$ and their punctual modifications. These arguments are rather standard.

\begin{prop}[Punctual nature] \label{local.nature.prop}
  Let $X$ be a Noetherian scheme, $Z\subset  X$ a closed subscheme and  ${\mathbf P}$  any of the 5 properties in
Definitions \ref{5.props.den}--\ref{5.props.den.1}. 
\begin{enumerate}
\item  If ${\mathbf P}$ holds for $Z\subset X$ then 
 ${\mathbf P}$ holds for  $(x, X_x)$ for every $x\in Z$.
\item The converse holds except possibly for TN. 
\end{enumerate}
\end{prop}

\noindent{\bf Complement.} We prove in (\ref{TN=STN.factor.cor}) that TN=STN, hence  (\ref{local.nature.prop}.2) also holds for TN.
\medskip

Proof.  Let
$f:Y\to X$ be a morphism that shows that $Z\subset X$ 
does not satisfy ${\mathbf P}$  and let $x$ be a generic point of
$\supp \bigl(f_*\o_Y/\o_X\bigr)$. Then localizing at $x$ shows that 
 $(x, X_x)$ does not satisfiy ${\mathbf P}$  in the cases N, SN and WN.
For STN, let  $g: Z\to X$ be the purely inseparable closure  of $X$ in $Y$ 
(\ref{weak.norm.2.defn}).  Then we take  
 $x$ to be  a generic point of
$\supp \bigl(f_*\o_Y/g_*\o_Z\bigr)$.

Let $f_x:(y,Y)\to (x,X)$ show that  $(x, X_x)$
does not satisfy ${\mathbf P}$. By Corollary \ref{relnorm.sheaf.cor},
it is obtained as the localization of a finite modification 
$f:Y\to X$ centered at
$\bar x\subset Z$. 
 If ${\mathbf P}$
is one of N, TN then  $f$ 
 shows that $Z\subset X$ 
does not satisfy ${\mathbf P}$.

If ${\mathbf P}$ is WN (resp.\ SN) then we have to find such an $f:Y\to X$ 
that is a universal homeomorphism (resp.\ also preserves residue fields).
These are also guaranteed by (\ref{relnorm.sheaf.cor}).
For STN we first extend $f_x$ to a finite morphism $f_U:Y_U\to U$  for
 some open neighborhood
 $x\in U\subset X$  and then use (\ref{extension.prop}.2) 
to extend it to $f:Y\to X$. 
\qed

\medskip

Even stronger localization can be obtained using the
following result which shows that the properties descend
for faithfully flat morphisms.

\begin{prop}\label{hart.ff.up.lem}
 Let $f:Y\to X$ be a faithfully flat morphism of Noetherian schemes, 
$Z\subset X$ a closed, nowhere dense subscheme and
$Z_Y:=f^{-1}(Z)$.  If $Z_Y\subset Y$ 
 satisfies any of the 5 properties ${\mathbf P}$ in
Definitions \ref{5.props.den}--\ref{5.props.den.1} then so does  $Z\subset X$. 
\end{prop}

Proof. Let $p:X'\to X$ be a finite modification centered at  $Z$. 
Then $p_Y:=(p\times f):Y':=X'\times_XY\to Y$ is a finite morphism that is an
isomorphism over $Y\setminus Z_Y$. 
Since $f$ is flat, none of the associated primes of
$Y'$ is contained in $p_Y^{-1}(Z')$. 
Thus  $p_Y:Y'\to Y$ is a finite modification
centered at $Z_Y$. Furthermore, if $p$ is a
universal homeomorphism  or  preserves residue fields
then $p_Y$ also has these properties.
If $f$ is faithfully flat then $p$ is an isomorphism off $p_Y$ is.

Thus if $p:X'\to X$ shows that $Z\subset X$ does not 
satisfy the property ${\mathbf P}$ then
$p_Y:Y'\to Y$ shows that $Z_Y\subset Y$ also does not 
satisfy the property ${\mathbf P}$.
 \qed

\begin{exmp} If $f:(y,Y)\to (x, X)$ is flat and $(y, Y)$  satisfies
the property ${\mathbf P}$ then $(x, X)$ need not satisfy ${\mathbf P}$.
For instance, if $(x,X)$ is any local scheme then the pair
$\{(x,0,0)\}\subset X\times \a^2$  is punctually normal.
\end{exmp}

The following could have been an exercise in \cite[Sec.II.5]{hartsh}.

\begin{prop}[Extension of finite morphisms]\label{extension.prop}
 Let $X$ be a Noetherian scheme, $X^0\subset X$ an open subscheme and $f^0:Y^0\to X^0$ a finite morphism. Then 
\begin{enumerate}
\item $f^0$ can be extended 
 to a finite, surjective morphism $f:Y\to X$.
\item Let $U\subset X$ be an open subset such that
 $f^0$ is an isomorphism (resp.\  partial weak or  seminormalization (\ref{weak.norm.defn})) over $U\cap X^0$.  
Then we can choose $f$ to be an isomorphism (resp.\  partial weak or  seminormalization) over $U$.
\item If $X$ is Nagata, $Y^0$  is reduced and  $f^0(Y^0)$
is dense in $X$  then there is a unique maximal reduced extension.
\end{enumerate}
\end{prop}

Proof.  By adding  $X\setminus \overline{X^0}$ both to $X^0$ and $Y^0$
we may assume that $X^0$ is dense in $X$. Then we would like to take
 $Y$ to be the
relative normalization of $X$ in $Y^0$; see (\ref{weak.norm.2.defn})
or \cite[Tag 0BAK]{stacks-project}. In case (3) this gives a finite
extension $Y\to X$ but in general the resulting  $\pi:Y^{\rm max}\to X$ is 
only integral. 
(For example, the relative normalization of  $ k[x,x^{-1}, y]/(y^2)$
in $ k[x,y]/(y^2)$ 
is $ k[x,x^{-r}y:r=0,1,\dots]/(y^2)$.)

Thus one needs to prove that a suitable finite subalgebra of
$\pi_*\o_Y^{\rm max}$ works;  the proof is very similar to 
\cite[Exrc.II.5.15]{hartsh} or see \cite[Tag 05K0]{stacks-project}.

If $f^0$ is an isomorphism over $U\cap X^0$ then first we extend it
to $U\cup Y^0$ (glued along $U\cap X^0$) and then apply (1). 
Finally,  if $f^0$ is a partial weak (resp.\ semi) normalization then first we choose 
any finite extension $f':Y'\to X$ and then let $Y$ be the relative  weak (resp.\ semi) normalization
(\ref{weak.norm.defn}) of $X$ in $Y'$. \qed
\medskip

If $Z'\subset X'$ is a localization of $Z\subset X$ and
$Y'\to X'$ is a finite modification then we can first extend
it to $Y^0\to X^0$ for some open subset $X^0\subset X$
and then use (\ref{extension.prop}) to extend it to $Y\to X$.
Thus (\ref{extension.prop}) can be  reformulated as follows.

\begin{cor} \label{relnorm.sheaf.cor} 
 The relative normalization
$X_Z^{\rm rn}$ commutes with localization.  The same holds for the
 relative  weak and seminormalization 
  $X_Z^{\rm rwn}$ and $X_Z^{\rm rsn}$ (\ref{weak.norm.2.defn}).
 \qed
\end{cor}

\section{Formal nature}

A key advantage of the punctual versions of the  4 properties in
Definition \ref{5.props.den} is that they are preserved by completions.
This is in contrast with the usual notions of normality
(resp.\ weak or seminormality); these are not always preserved by completions
and  it is frequently not easy to understand the cases when they are preserved.

The following  is a special case of formal gluing originated in  \cite{artin}.
For a thorough general discussion see 
\cite[Tags 05E5 and 0AEP]{stacks-project}. 
We give a short proof below  using only the simpler methods employed there. 

\begin{prop} \label{sections.completions.cor}
 Let $X$ be a  Noetherian scheme and $Z\subset X$ a
closed, nowhere dense  subscheme with ideal sheaf $I\subset \o_X$.
Let $\hat X$ denote the $I$-adic completion of $X$ and
$Z\cong \hat Z\subset \hat X$ the corresponding subscheme.
Then completion provides an equivalence of categories
$$
\{\mbox{finite modifications of $Z\subset X$}\}
\Leftrightarrow
\{\mbox{finite modifications of $\hat Z\subset \hat X$}\}.
$$
\end{prop}

Proof. The question is local hence we may assume that $X$ is affine.
We may also assume that none of the associated primes of $X$ is contained in
$Z$. 

Let $\tau:U:=X\setminus Z\into X$ be the natural open embedding;
similarly we have $\hat\tau:\hat U:=\hat X\setminus \hat Z\into \hat X$.
Set $R:=H^0(X, \o_X), S:=H^0(U, \o_U)=H^0(X, \tau_*\o_U)$ and, similarly
$\hat R:=H^0(\hat X, \o_{\hat X}), \hat S:=H^0(\hat U, \o_{\hat U})=
H^0(\hat X, \hat\tau_*\o_{\hat U})$.

Since $\hat X\to X$ is flat, cohomology and base change 
for $\tau:U\to X$ says that
$\hat S=S\otimes_R\hat R$.
 
As we noted in (\ref{fin.mod.defn}.2), finite modifications of $Z\subset X$
correspond to subalgebras $R\subset R'\subset S$ that are finite over $R$.

In general, let $R$ be a  Noetherian  ring, $I\subset R$  an ideal and $S\supset R$ an $R$-algebra
such that $I$ is nilpotent on $S/R$. If $R\subset R'\subset S$
is a  subalgebra that is finite over $R$ then  $I^mR'\subset R$ 
for some $m>0$ and hence $R'\subset R:I^m$.
Multiplication gives a map
$$
\phi: (R:I^m)/I^m\otimes (R:I^m)/I^m\to (R:I^{2m})/R
\qtq{and}
\phi\bigl(R'/I^m, R'/I^m\bigr)\subset R'/R.
$$
Conversely, if $R\subset M\subset (R:I^m)$ is any submodule such that
$\phi\bigl(M/I^m, M/I^m\bigr)\subset M/R$ then $M$ is a subalgebra.

Take now the $I$-adic completion $\hat R$ and set $\hat S:=\hat R\otimes_RS$.
Since $I$ is nilpotent on $S/R$, we see that $\hat S/\hat R\cong S/R$
and so $(R:I^m)/I^m\cong (\hat R:\hat I^m)/\hat I^m$ for every
 $m>0$.
This gives an equivalence of categories between
algebras $R\subset R'\subset S$ that are finite over $R$
and algebras $\hat R\subset \hat R'\subset \hat S$ that are finite over $\hat R$. \qed

\medskip

By passing to the limit on both sides, we get the following.

\begin{cor} Let $X$ be a  Noetherian scheme and $Z\subset X$ a
closed, nowhere dense  subscheme. Then
$$
\hat X_{\hat Z}^{\rm rn}=X_Z^{\rm rn}\times_X\hat X. \qed
$$
\end{cor}

\begin{cor}[Formal nature]\label{formal.nature.thm}
Using the notation of (\ref{sections.completions.cor}) 
let   ${\mathbf P}$  denote any of the 4 properties in
Definition \ref{5.props.den}.  Then ${\mathbf P}$ holds for $Z\subset X$ iff 
 it holds for  $\hat Z\subset \hat X$.
\end{cor}

\noindent{\bf Complement.} We prove in (\ref{TN=STN.factor.cor}) that STN=TN, hence (\ref{formal.nature.thm})  also  holds
for STN.
\medskip

Proof.  The if part follows from (\ref{hart.ff.up.lem}) since
$\hat X$ is faithfully flat over $X$.

To see the converse, 
let $\hat p:\hat Y\to \hat X$ be any finite modification. By
(\ref{sections.completions.cor})   there is a finite modification $p: Y\to X$ such that
$\hat Y=Y\times_X\hat X$. Since $\hat X$ is faithfully flat over $X$,
$\hat p$ is an isomorphism (resp.\ universal homeomorphism) iff $p$ is.
\qed
\medskip

The above result, together with (\ref{pf.of.hrt.defn.2}) and (\ref{local.nature.prop}) gives the following.

\begin{cor} \label{n.w.s.punctual.cor}
A Noetherian  scheme $X$ is normal (resp.\ weakly or seminormal) iff the completion  $(\hat x, \hat X)$ is punctually normal (resp.\ weakly or seminormal) for every $x\in X$. \qed
\end{cor}

\section{Finiteness properties}

Next we study when a relative normalization is finite over $X$.

\begin{thm}\label{punc.nag.lem}
Let $X$ be a  Noetherian  scheme and  $Z\subset X$ a closed, 
nowhere dense subscheme. The following are equivalent.
\begin{enumerate}
\item The relative normalization $\pi_Z: X_Z^{\rm rn}\to X$ is finite. 
\item The punctual normalization 
$\pi_x: \bigl(x^{\rm pn},X^{\rm pn}\bigr)\to (x,X)$ is finite for every  $x\in  Z$. 
\end{enumerate}
\end{thm}

Proof.  (1) $\Rightarrow$ (2)  since 
$\pi_x: \bigl(x^{\rm pn},X^{\rm pn}\bigr)\to (x,X)$ is the localization of
$\pi_{\bar x}: X_{\bar x}^{\rm rn}\to X$  by (\ref{relnorm.sheaf.cor})
 and $X_Z^{\rm rn}$ dominates
 $X_{\bar x}^{\rm rn}$ if  ${\bar x}\subset Z$.

Conversely, assume that (2) holds. We may assume that none of the associated 
points of $X$ are contained in $Z$.

We repeatedly use the following construction.
For  any point $x\in X$ take the 
punctual normalization  $\bigl(x^{\rm pn}, X^{\rm pn}\bigr)\to (x, X)$
and then extend it to a finite modification  $p_1:X_1\to X$ centered at
$\bar x$ using
(\ref{extension.prop}). As we note in (\ref{PN1.updown.lem}),
the assumptions of (2) also hold for $X_1$ and $Z_1:=p_1^{-1}(Z)$.

The canonical choice would seem to be to apply  this construction first to the
closed points where $X$ is not punctually normal. This works well
and we get $X_1$ that is  punctually normal at all 
closed points. However, when we apply the procedure to $X_1$,
we eliminate the 1-dimensional points where $X_1$ is not 
 punctually normal but we may generate new  punctually non-normal closed points. Thus we have to proceed from the other end. The price we pay is that
the procedure itself is non-canonical. 

First  we apply the construction  to the   points of $Z$ that have codimension 1 in $X$. 
We get $p_1:X_1\to X$ 
such that $X_1$ is regular at all  
 points of $Z_1$ that have codimension 1 in $X_1$. By 
(\ref{nonrom.finite.inZ.lem}) there are only finitely   points $x_1\in Z_1$ such that
the pair $(x_1, X_1)$ is not normal.
The closures of these points form a 
closed subset $W_1$. Next we repeat the procedure for the 
generic points of $W_1$ to get $p_2:X_2\to X_1$. 
After finitely many steps
we get  a finite partial modification $\pi:X_r\to X$ 
centered at $Z$ 
such that  $X_r$ is 
punctually normal at all points of $\pi^{-1}(Z)$. Thus 
 $\pi:X_r\to X$ is the relative normalization of the pair
 $Z\subset X$ by Proposition \ref{local.nature.prop}.
\qed

\begin{lem} \label{nonrom.finite.inZ.lem}
Let $X$ be a Noetherian  scheme and $W\subset X$ a closed
subscheme that does not contain any of the associated points of $X$.
Then there are only  finitely many points $x\in W$ such that
the pair $(x, X)$ is not normal.
\end{lem}

Proof. The question is local so we may assume that $X$ is affine.
By our assumptions there is  a Cartier divisor $(g=0)$   containing $W$. 
 By Lemma \ref{hart.serre.loc.thm}, if a pair $(x, X)$ is not normal
then either $x$ has codimension 1 in $X$  (thus $x$ is a generic point of $Z$)
or  $\depth_xX=1$ and $\codim_Xx\geq 2$. By
Corollary \ref{cart.pair.nomal.cor} such an  $x$ is an associated point of $(g=0)$. Since $X$ is Noetherian,
there are only finitely many such points.\qed

\begin{defn} \label{PN1.defn}
A  Noetherian  scheme $X$ is called
{\it punctually N-1} at $x\in X$ if the punctual normalization 
$\pi_x: \bigl(x^{\rm pn},X^{\rm pn}\bigr)\to (x,X)$ is finite.
$X$  is called
{\it punctually N-1} if this holds for every $x\in X$.

If $X$ is reduced then the normalization of a local ring dominates 
its punctual normalization, hence if the local rings  are
N-1 (\ref{nagata.defn}) then $X$ is also punctually N-1.
(The converse probably does not hold; see (\ref{N-1.PN-1.fin.prop}) and
(\ref{ques.1}).)
\end{defn}

The following is closely related to
\cite[Tag 0333]{stacks-project}.

\begin{prop} \label{N-1.PN-1.fin.prop}
A Noetherian integral scheme $X$ is N-1 iff it is
punctually N-1 and there are only finitely many points $x_i\in X$
such that $(x_i, X)$ is not punctually normal.
\end{prop}

Proof.  Let $X^{\rm n}\to X$ denote the  normalization
with structure sheaf $\o_{X}^{\rm n}$ and 
$X_i^{\rm pn}\to X$  the punctual normalization at $x_i$
with structure sheaf  $\o_{X_i}^{\rm pn}$. Since $X^{\rm n}$ dominates
$X_i^{\rm pn}$ we get injections
$$
\o_{X_i}^{\rm pn}/\o_X\into \o_{X}^{\rm n}/\o_X.
$$
If $\o_{X}^{\rm n}$ is a coherent $\o_X$-sheaf then 
$\o_{X_i}^{\rm pn}$ is also coherent and there are only finitey many of them since
each $x_i$ is an associated point of $\o_{X}^{\rm n}/\o_X $.

Conversely, assume that $X$ is punctually N-1 and there are only finitely many points $x_i\in X$
such that $(x_i, X)$ is not punctually normal. Apply
(\ref{punc.nag.lem}) with $Z:=\cup_i\bar x_i$. The relative normalization
$\pi_Z:X_Z^{\rm rn}\to X$ is finite and punctually normal along the preimage of $Z$. Since $\pi_Z$ is an isomorphism over $X\setminus Z$, we conclude that
$X_Z^{\rm rn} $ is punctually normal everyhere hence normal by Propositions \ref{pf.of.hrt.defn.2} and \ref{local.nature.prop}. \qed

\begin{lem}  \label{PN1.updown.lem}
Let  $g:Y\to X$ be a finite modification  of Noetherian   schemes
centered at some $W\subset X$. Assume that none of the associated points of
$X$ is contained in $W$.  Pick a point $x\in X$ and set $y:=\red g^{-1}(x)$.
Then  $X$ is punctually N-1 at $x$ iff  $Y$ is punctually N-1 at $y$.
\end{lem}

Proof. If $x\not\in W$ then  $X_x\cong Y_y$ so the claim is interesting
only if $x\in W$.

$Y^{\rm pn}$ dominates $X^{\rm pn}$, 
so the only if part is clear. 
To see the converse, we may assume that $X$ is affine. Then
 there is a non-zerodivisor  $\phi\in \o_{X}$ such that
$\phi\cdot \o_Y\subset \o_X$. Thus $\phi\cdot \o_Y^{\rm pn}\subset \o_X^{\rm pn}$
hence if $\o_X^{\rm pn}$ is coherent then so is $\o_Y^{\rm pn}$.\qed
\medskip

Next we study the effect of nilpotents on the relative normalization.

\begin{say}[Structure of the relative normalization]
\label{srt.of.relnorm.say}
Let $X$ be a Noetherian scheme,
$Z\subset  X$ a closed, nowhere dense  subscheme and $j:U:=X\setminus Z\into X$
the natural open embedding.  Assume for notational simplicity that 
none of the associated primes of $X$ is contained in $Z$. 
 Let $J\subset \o_X$ be a nilpotent ideal.
Note that $j_*(J|_U)\subset j_*\o_U$ is nilpotent, hence its sections are integral over $X$. Thus we have  exact sequences
$$
\begin{array}{ccccccccc}
0 &\to & j_*(J|_U)&\to &  j_*\o_U&\stackrel{\rho}{\to}&  j_*(\o_U/J|_U)
&\to& R^1j_*(J|_U)\\
0 &\to & j_*(J|_U)&\to &  \o_X^{\rm rn}&\stackrel{\rho}{\to}&  (\o_X/J)^{\rm rn} 
&\to& R^1j_*(J|_U)\\
\end{array}
\eqno{(\ref{srt.of.relnorm.say}.1)}
$$
where $\o_X^{\rm rn}$ and $(\o_X/J)^{\rm rn}$ denote the
structure sheaves of the relative normalizations of the pairs
$Z\subset  X$ and   $Z\subset  \spec_X(\o_X/J)$.
It is rather unclear when $\rho$ is surjective; see 
(\ref{tnr.not.top.exmp}) for an example. However, if  $R^1j_*(J|_U)=0$  then $\rho$ is surjective and we obtain the following.
\medskip

{\it Claim \ref{srt.of.relnorm.say}.2.} If $\dim\supp J=1$ then
every finite modification of  the pair $Z\subset \spec_X(\o_X/J)$ lifts to a
finite modification of  $Z\subset X$. 
In particular,   $Z\subset X$ satisfies TN  iff
$Z\subset \spec_X(\o_X/J)$ does. \qed
\medskip

The second sequence in  (\ref{srt.of.relnorm.say}.1) shows that 
 if $ \o_X^{\rm rn} $ is coherent then so is $j_*(J|_U)$. 
By the easy direction of  \cite[IV.5.11.1]{EGA}
(or by the argument in Proposition \ref{pf.of.hrt.defn.2})
this implies the following.

\medskip

{\it Claim \ref{srt.of.relnorm.say}.3.} Assume that  $ X_Z^{\rm rn}\to X $ is 
finite
and let $J\subset \o_X$ be a nilpotent ideal. Then either $\supp J\subset Z$ or $Z\cap\supp J$ has codimension $\geq 2$ in $\supp J$.

Thus if  $ X_Z^{\rm rn}\to X $ is 
finite for every $Z$ then  the nilradical of $\o_X$ has zero-dimensional support.
\qed
\medskip

Next let $I\subset \o_X$  be the nilradical
 and set
$\o_X^{\rm nil}:=\rho^{-1}\bigl(\o_{\red X}\bigr)$. 
 We have an exact sequence
$$
0\to j_*(I|_U)\to  \o_X^{\rm nil}\stackrel{\rho}{\to}  \o_{\red X}
\eqno{(\ref{srt.of.relnorm.say}.4)}
$$
and $\pi:\spec_X  \o_X^{\rm nil}\to X$ is a limit of 
residue field preserving universal
homeomorphisms. Note that $\pi$  is an isomorphism iff $I=j_*(I|_U)$; equivalently, if
$\depth_ZI\geq 2$. Thus we obtain the following.

\medskip

{\it Claim \ref{srt.of.relnorm.say}.5.}
Let $Z\subset X$ be a seminormal pair  and $I\subset \o_X$  the nilradical. Then $\depth_ZI\geq 2$ and 
so $\codim_W(Z\cap W)\geq 2$  for every  associated prime $W$   of $I$. 
\qed
\end{say}

The implication of (\ref{srt.of.relnorm.say}.3) becomes an
equivalence  for Nagata schemes  (\ref{nagata.defn}).

 \begin{prop} \label{relnorm.finite.crit.prop}
Let $X$ be a Nagata scheme  and
$Z\subset  X$ a closed, nowhere dense  subscheme. The following are equivalent.
\begin{enumerate}
\item The relative normalization $X_Z^{\rm rn}\to X$ is finite over $X$.
\item Let $J\subset \o_X$ be a nilpotent ideal. Then either $\supp J\subset Z$ or $Z\cap\supp J$ has codimension $\geq 2$ in $\supp J$.
\end{enumerate}
\end{prop}

Proof.  (1) $\Rightarrow$ (2) follows from (\ref{srt.of.relnorm.say}.3).

To see the converse, let $I\subset \o_X$ be the  nilradical.
By assumption, if $W\subset X$ is any associated prime of $I$ then
either $W\subset Z$ or  $Z\cap W$ has codimension $\geq 2$ in $W$. Thus $j_*(I|_U)$ is
coherent by \cite[IV.5.11.1]{EGA}. The normalization $(\red X)^{\rm n}$ is finite over
$\red X$ since $X$ is Nagata and (\ref{srt.of.relnorm.say}.1)
gives the exact sequence
$$
0\to j_*(I|_U)\to  \o_X^{\rm rn}\stackrel{\rho}{\to}  \o_{\red X}^{\rm n} . \qed
$$

 A complete semilocal ring is  Nagata  \cite[32.1]{MR0155856}
and,  again using   \cite[IV.5.11.1]{EGA}, we see that $j_*\o_U$ is coherent 
iff $ X$ does not have 
 1-dimensional irreducible components. 
 Hence we have
 proved the following.

\begin{cor} \label{relnorm.finite.crit.cor}
 Let $(x,X)$ be a Noetherian, semilocal, complete scheme
without isolated points
and $j:U:=X\setminus\{x\}\into X$ the natural embedding.
\begin{enumerate}
\item The punctual normalization
 $\pi:\bigl(x^{\rm pn}, X^{\rm pn}\bigr)\to (x, X)$ is finite 
iff $\o_X$ has no nilpotent ideals with 1-dimensional support. 
\item The punctual hull $j_*\o_U$ is integral over $\o_X$ 
iff $ X$ does not have 
 1-dimensional irreducible components.
\item The punctual hull $j_*\o_U$ is finite over $\o_X$
 iff $ X$ does not have 
 1-dimensional associated primes.\qed
\end{enumerate}
\end{cor}

Combining  (\ref{relnorm.finite.crit.cor}) and  (\ref{sections.completions.cor})
 gives a characterization of punctually N-1 local schemes.
The 1-dimensional case is essentially in \cite{MR1512631};
see also \cite[Sec.1.13]{k-res}.

\begin{cor} \label{relnorm.finite.crit.cor.2}
 Let $(x,X)$ be a Noetherian, semilocal scheme without isolated points
and $j:U:=X\setminus\{x\}\into X$ the natural embedding.
\begin{enumerate}
\item The punctual normalization
 $\pi:\bigl(x^{\rm pn}, X^{\rm pn}\bigr)\to (x, X)$ is finite 
iff the completion $\hat X$ has no nilpotent ideals with 1-dimensional support. 
\item The punctual hull $j_*\o_U$ is integral over $\o_X$ 
iff $\hat X$ does not have 
 1-dimensional irreducible components.
\item The punctual hull $j_*\o_U$ is finite over $\o_X$
 iff $\hat X$ does not have 
 1-dimensional associated primes.\qed
\end{enumerate}
\end{cor}

These give a  characterization of  weakly and seminormal schemes in terms of the
conductor of the normalization. This is well known if the normalization is finite (\ref{wn.sn.conductor.crit}), but in general 
 the (global) conductor ideal could be trivial.
We go around this problem by working punctually.

\begin{cor} \label{sn.wn.loc.char.cor} For a Noetherian scheme $X$ the following are equivalent.
\begin{enumerate}
\item  $X$ is semi  (resp.\ weakly) normal.
\item $X$  is reduced and for every $x\in X$, the punctual normalization
$\pi_x:X_x^{\rm pn}\to X$ is finite and its conductor 
$C_x^{\rm pn}\subset X_x^{\rm pn}$ is reduced (resp.\ $\{x\}$ is 
purely inseparably closed in $C_x^{\rm pn}$).
\end{enumerate}
\end{cor}

 Note that in (2) the   conductor 
$C_x^{\rm pn}$ is 0-dimensional, so these are very simple conditions.
\medskip

Proof. If $X$ is seminormal then the punctual normalization  
 $(x^{\rm pn},X^{\rm pn})\to (x,X)$ is finite by
(\ref{srt.of.relnorm.say}.5), (\ref{formal.nature.thm}) and (\ref{relnorm.finite.crit.cor.2}).
(This  should also follow from a theorem of Chevalley and Mori; see \cite[33.10]{MR0155856}.) Then (2) holds 
by (\ref{wn.sn.conductor.crit}). 

Conversely, if (2) holds then $(x, X)$ satisfies SN (resp.\ WN) for  every $x\in X$ by (\ref{wn.sn.conductor.crit})
hence $X$ is semi  (resp.\ weakly) normal by (\ref{local.nature.prop}). \qed

\medskip

Observe that if $X$ is a Noetherian, integral, seminormal  scheme whose normalization is not finite then, by (\ref{N-1.PN-1.fin.prop}),  there are infinitely many distinct points  $x_i\in X$ such that $(x_i, X)$ is not
punctually normal.

\section{Factoring topological modifications}

\begin{thm}  \label{factor.topmod.thm}
Let $X$ be a  Noetherian scheme, $Z\subset X$ a
closed, nowhere dense  subscheme and  $g:Y\to X$ a finite,
topological modification centered at $Z$. Then $g$ can be factored as
$$
g:Y\stackrel{g'}{\to} X'\stackrel{\pi}{\to} X
$$
where $g'$ is a finite, universal homeomorphism and
$\pi$ is a finite  modification centered at $Z$.
\end{thm} 

\begin{rem} The factorization is not unique and there does not seem to be 
a canonical choice for it, except when $\dim Z=0$  (\ref{ptn.say}). 
For example,  $\c[x^2,  xy, y]\to \c[x, y]$ can be factored as
$$
\c[x^2,  xy, y] \to  \c[x^2, x^{2n+1}, xy, y]\to \c[x, y]
$$
for any $n\geq 0$, thus there is no minimal choice for $X'$. 
Similarly, $\c[x, y]/(x^2y^2)\to \c[x]+\c[y]$ can be factored as
$$
\c[x, y]/(x^2y^2)\to \c[x, x^{-n}\epsilon]+\c[y,y^{-n}\epsilon]\to \c[x]+\c[y]
$$
for any $n\geq 0$, thus there is no maximal choice for $X'$. 
\end{rem}

\begin{cor} \label{TN=STN.factor.cor}
 Let $X$ be a  Noetherian scheme and $Z\subset X$ a
closed, nowhere dense  subscheme. Then $Z\subset X$  satisfies TN iff it
satisfies STN. 
\end{cor}

Proof. It is clear that STN $\Rightarrow$ TN.
Next assume that we have $g:Y\to X$ showing that STN fails for $Z\subset X$.
Using (\ref{factor.topmod.thm}) we get $\pi:X'\to X $ which shows that  TN fails for 
$Z\subset X$. \qed

\begin{cor} \label{TN=unibranch.cor}
For a Noetherian scheme $X$ the  following are equivalent.
\begin{enumerate}
\item $X$ is  topologically normal,
\item $X$ is geometrially unibranch (see \cite[Sec.IV.6.15]{EGA} or \cite[Tag 06DJ]{stacks-project}),
\item the normalization $X^{\rm n}\to X$  is a universal homeomorphism,
\item the weak normalization  $X^{\rm wn}$ is normal.\qed
\end{enumerate}
\end{cor}

\begin{say}[Punctual topological normalization]\label{ptn.say}

Let $(x, X)$ be a complete local scheme and 
 $p:(y, Y)\to (x, X)$  a finite morphism. 
Let $A\subset k(y)$ be the separable closure of $k(x)$. 
There is a (unique)  finite, \'etale morphism  $ (x', X')\to (x, X)$ 
whose fiber over $x$ is isomorphic to $\spec A$. Furthermore,
$p$ lifts to $p_Y:(y, Y)\to (x', X')$ since $ (x', X')\to (x, X)$  is formally smooth and $(y, Y)$ is complete.

Assume next that $p:Y\setminus\{y\}\to X\setminus\{x\}$ is a
universal homeomorphism.  Then $\red(p_Y(Y))$ is a union of
some of the irreducible components of $\red(X')$. 
Thus there is a largest subscheme $X^u\subset X'$ such that
$\red(X^u)=\red(p_Y(Y))$ and $x$ is not an associated point of  $X^u$. 
Then $X^u\setminus\{x'\}\to X\setminus\{x\}$ is \'etale and
an isomorphism on the underlying reduced subschemes, hence an isomorphism.

We have thus factored $p$ as
$$
p: (y, Y)\stackrel{p'}{\to} (x^u, X^u)\stackrel{p^u}{\to}  (x, X)
\eqno{(\ref{ptn.say}.1)}
$$
where $p^u$ is an unramified punctual modification and  $p'$ 
is a universal homeomorphism. 
By (\ref{sections.completions.cor}), this factorization also exists  even if 
$(x, X)$ is  not  complete. This proves (\ref{factor.topmod.thm})
in the special case when  $Z$ is a closed point.

Applying this to  
$p:\bigl(\hat x^{\rm pn}, (\red\hat X)^{\rm pn}\bigr)\to (\hat x, \red\hat X)\to (\hat x, \hat X)$ and then descending to $X$ we obtain the following.
\medskip

{\it Claim \ref{ptn.say}.2.}  Let $(x, X)$ be a semilocal, Noetherian scheme.
Then there is a unique punctual modification
$$
\pi^{\rm ptn}: \bigl(x^{\rm ptn}, X^{\rm ptn}\bigr)\to (x, X)
$$
such that $\pi^{\rm ptn}$ is unramified and $ \bigl(x^{\rm ptn}, X^{\rm ptn}\bigr)$
is  the smallest 
topologically normal punctual modification of $(x, X)$. \qed
\medskip

Note that if $(x, X)$ is Henselian with perfect residue field then
$X^{\rm ptn}$ is simply the disjoint union of the closures of the connected  components of
$X\setminus\{x\}$; cf.\  (\ref{hesnelian.topext.conn.prop}). De~Jong pointed out that  (\ref{ptn.say}.2) can also be obtained from this by Galois descent.

\end{say}

\begin{say}[Proof of Theorem \ref{factor.topmod.thm}]
\label{factor.topmod.thm.pf}
We use Noetherian induction on $g^{-1}(Z)\subset Y$. 

Let $x\in X$ be the generic points of $Z$. After localizing at $x$ we obtain
  a punctual, topological modification $p:(y, Y)\to (x, X)$.
By (\ref{ptn.say}.1) it can be factored  as
$$
p: (y, Y)\stackrel{p'}{\to} (x', X')\stackrel{p''}{\to}  (x, X),
\eqno{(\ref{factor.topmod.thm.pf}.1)}
$$
where $p''$ is an unramified punctual modification and  $p'$ 
is a universal homeomorphism.  We extend (\ref{factor.topmod.thm.pf}.1)
first to an open neighborhood and then, using  (\ref{extension.prop}),
we get morphisms
$$
g:Y\stackrel{g_1}{\to} X_1\stackrel{p_1}{\to}X,
\eqno{(\ref{factor.topmod.thm.pf}.2)}
$$
where $p_1$ is a finite modification centered at $Z$ and
$g_1:Y\to X_1$ is a finite morphism that is a
universal homeomorphism outside $p_1^{-1}(Z)$ and also
over the generic points of $p_1^{-1}(Z)$.
By Noetherian induction,   $g_1$ can be factored as
$$
g_1:Y\stackrel{g'_1}{\to} X'_1\stackrel{\pi_1}{\to} X_1
\eqno{(\ref{factor.topmod.thm.pf}.3)}
$$
where $g'_1$ is a finite, universal homeomorphism and
$\pi_1$ is a finite  modification centered at $p^{-1}(Z)$.
Thus we can set $X':=X'_1$, $g':=g'_1$   and 
$$
g:Y\stackrel{g'}{\to} X'\stackrel{p_1\circ \pi_1}{\longrightarrow} X
\eqno{(\ref{factor.topmod.thm.pf}.4)}
$$
gives the required factorization for $g$. \qed

\end{say}

Arguing as in 
(\ref{punc.nag.lem}) and using (\ref{ptn.say}.2), this implies the following.

\begin{cor}  Let $X$ be a  Noetherian scheme and $Z\subset X$ a
closed, nowhere dense  subscheme. Then there is a    finite
 modification $\pi:Y\to X$ centered at $Z$ such that the pair
$\pi^{-1}(Z)\subset Y$ is topologically normal. \qed
\end{cor}

The following example shows that there does not seem to be a good notion of
strong topological normality with residue field preservation.

\begin{exmp} \label{tnr.not.top.exmp}
Start with the function field $K=k(t)$ with $\chr k=2$. Set
$$
R:=K[x,y,\sqrt{t}x, \sqrt{t}y]\subset K(\sqrt{t})[x,y].
$$
Note that the normalization of $R$ is $K(\sqrt{t})[x,y]$ and 
the conductor ideal is  $(x,y)\subset K(\sqrt{t})[x,y]$. Thus $R$ is 
seminormal but not weakly normal and does not satisfy TNR.

We can embed $\spec R$ into $\a^4_K$ by $(x,y,\sqrt{t}x, \sqrt{t}y)
\mapsto (u_1, v_1, u_2, v_2)$. The image is defined by the obvious equations
$$
u_1v_2-u_2v_1= tu_1v_1-u_2v_2= tu_1^2-u_2^2= tv_1^2-v_2^2=0.
$$
Let $X\subset \a^4_K$ be defined by the last 2 equations. 
It is a complete intersection, hence $S_2$ and so
$\{\mathbf{0}\}\subset X$ is a normal pair which also satisfies TNR. 

On the other hand, the last 2 equations imply that
$$
(u_1v_2-u_2v_1)^2= (tu_1v_1-u_2v_2)^2=0,
$$
hence $\red X\cong \spec R$. 
Thus
$\{\mathbf{0}\}\subset X$  satisfies TNR but
$\{\mathbf{0}\}\subset \red X$ does not.
\end{exmp}

\section{Flat families}

All 5 properties behave  well for flat families.

\begin{thm} \label{reduced.morph.preserves.P}
Let $f:Y\to X$ be a flat morphism of Noetherian schemes
with (geometrically) reduced fibers
and  ${\mathbf P}$  any of the 5 properties in
Definitions \ref{5.props.den}--\ref{5.props.den.1}.
Assume that  $X$ and the geometric generic fibers satisfy
 ${\mathbf P}$. Then $Y$ also satisfies ${\mathbf P}$. 
\end{thm}

For normality this is classical; see for instance   \cite[23.9]{mats-cr}.
For seminormality this is proved in  \cite{MR571055}, for
weak normality in    \cite{MR556308}; in both cases for     N-1 schemes  
(called Mori schemes in these papers).
The argument below uses the method of   \cite{MR571055}.
\medskip

Proof.
Pick any point $y\in Y$ and set $x:=f(y)$. If $Y$ is contained in a generic fiber then
$(y,Y)$ satisfies ${\mathbf P}$ by assumption. If $y$ is a
generic point of its fiber $Y_x$ then, as we prove in (\ref{reg.morph.preserves.P}),  $(y,Y)$ satisfies ${\mathbf P}$. 

 $X$ is reduced by Proposition \ref{pf.of.hrt.defn.2}  in cases N, SN, WN and 
for every other point we have
$$
\depth_yY=\depth_xX+\depth_yY_x\geq 1+1=2,
\eqno{(\ref{reduced.morph.preserves.P}.1)}
$$
hence $(y,Y)$ is even punctually normal by Lemma \ref{hart.serre.loc.thm}
and it satisfies ${\mathbf P}$ for every ${\mathbf P}$.

As we noted in (\ref{basic.def.comms}.4), in case TN 
it is enough to show that $\red Y$ satisfies TN. 
 We may replace $X$ and $Y$ by $\red X$ and $\red Y$.  We do not actually need that
$\red X$  satisfies TN. It is $S_1$ and this is all we used
in (\ref{reduced.morph.preserves.P}.1).
\qed

\begin{prop}\label{reg.morph.preserves.P}
 Let $f:Y\to X$ be a flat, regular morphism of Noetherian schemes.
If $X$ satisfies any of the 5 properties ${\mathbf P}$ in
Definitions \ref{5.props.den}--\ref{5.props.den.1} then so does $Y$.
\end{prop}

Proof. Pick a point $y\in Y$ and set $x:=f(y)$. By (\ref{formal.nature.thm}),
$\hat X_x$ satisfies  ${\mathbf P}$. Note that
$Y\times_X\hat X_x\to Y$ is faithfully flat. Thus if
$\bigl(y, Y\times_X\hat X_x\bigr)$ satisfies ${\mathbf P}$
then so does $(y, Y)$ by Proposition \ref{hart.ff.up.lem}. Thus it is
sufficient to prove (\ref{reg.morph.preserves.P}) in case
$(x, X)$ is local and complete. 
We need slightly different arguments in the various cases.

If $(x, X)$ is punctually normal then either $X$ is regular
(of dimension 1) hence $Y$ is also regular, or
$\depth_xX\geq 2$ hence $\depth_yY\geq 2$.

For TN  we may assume, using  (\ref{srt.of.relnorm.say}.2),  that $\o_X$ has no nilpotent ideals with 1-dimensional support, and then 
the punctual normalization  
 $(x^{\rm pn},X^{\rm pn})\to (x,X)$ is finite by (\ref{relnorm.finite.crit.cor}).
It  is also a  universal homeomorphism  by our assumptions. Thus
$Y\times_XX^{\rm pn}\to Y$ is also a  universal homeomorphism  (resp.\ also preserves residue fields). By the already established normal case 
the pair
$ Y\times_Xx^{\rm pn}\subset Y\times_X X^{\rm pn}$ is normal.  Therefore,
as we noted in Paragraph \ref{props.from.relmorm.say.1}, $Y$
satisfies  TN.

For SN and WN we use
(\ref{sn.wn.loc.char.cor}),   thus the punctual normalization
$\pi:X^{\rm pn}\to X$ is finite and its conductor 
$C_X^{\rm pn}\subset X^{\rm pn}$ is reduced (resp.\ $\{x\}$ is 
purely inseparably closed in $C_X^{\rm pn}$). Therefore the conductor 
 $C_Y^{\rm pn}$
of $Y\times_X X^{\rm pn}\to Y$ is also reduced and  $Y_x:=Y\times_Xx$ is 
purely inseparably closed in it in the weakly normal case by (\ref{weak.closure.field case.exmp}). Thus $Y$ satisfies WN (resp.\ SN) by (\ref{wn.sn.conductor.crit}).\qed

\begin{exmp}\label{top.non-flat.fibers.say}
It is natural to ask what happens if in (\ref{reduced.morph.preserves.P})
we only assume that the fibers are
topologically normal in codimension $m$ for some $m$.
There are easy counter examples.

The fibers of $\spec k[x,y]/(y^2-x^2)\to \spec k[x]$ are 0-dimensional,
hence topologically normal yet $\spec k[x,y]/(y^2-x^2)$ is not
topologically normal.

Similarly, let $R\subset k[x,y]$ consist of those polynomials for which
$p(0,0)=p(0,1)$. Then the fibers of $\spec R\to \spec k[x]$ are 1-dimensional
hence topologically  normal in codimension $2$ yet $\spec R$ is not
topologically normal in codimension $2$.

Nonetheless, these are essentially the only counter examples.
We can even keep track of the codimension as in
(\ref{top.norm.m.thm}).
\end{exmp}

\begin{prop} Fix $n, m\geq 1$. Let $g:X\to S$ be a morphism of pure relative dimension $n$ 
for some $n$ and $W\subset X$ a closed subset.
Assume that 
\begin{enumerate}
\item $X\setminus W$ is topologically normal in codimension $m$.
\item every fiber $X_s$ is topologically normal  in codimension $m$ and
\item $\codim_{X_s}(W\cap X_s)\geq m$  for every $s\in S$.
\end{enumerate}
Then $X$ is topologically normal  in codimension $m$.
\end{prop}

Proof. Let  $f:Y\to X$ be a finite modification centered at $Z\subset X$
that is a putative counter example.  We need to show that $f$ is a 
 universal homeomorphism. By assumption (1) this holds over
$X\setminus W$, we can thus assume that $Z\subset W$.

We claim that $Y_s$ has pure dimension $n$ for every $s\in S$.
  It is enough to prove the claim
when $S$ is local and irreducible. If  $s\in S$ has codimension  $r$
then $s$ is set-theoretically the complete intersection of $r$
Cartier divisors, thus $Y_s$ is also 
set-theoretically the complete intersection of $r$
Cartier divisors, hence every irreducible component of $Y_s$ 
has dimension $\geq n$.
Since $X_s$ has pure dimension $n$ and $Y_s\to X_s$ is finite,
this implies that $Y_s$ has pure dimension $n$.

In particular,  none of the irreducible components of $Y_s$
is contained in $f_s^{-1}(W\cap X_s)$ hence
$ f_s:Y_s\to X_s$ is a finite topological modification centered at $W\cap X_s$.
Then $ f_s$ is a universal homeomorphism by assumption (2)
and so is $f$.\qed

\section{Connectedness properties}

A version of topological normality first appeared in \cite{MR0142547} where it is proved that $\depth_xX\geq 2$ implies that $(x, X)$ is punctually
TN. More generally, one can understand topological normality in
codimension $\geq m$ in terms of connectedness properties of \'etale
covers of $X$. Many parts of the following theorem are discussed in
\cite{MR0142547} and \cite[III.3 and XIII.2]{sga2}. 

\begin{thm}\label{top.norm.m.thm} Fix a natural number $m\geq 1$.
 For a Noetherian  scheme  $X$ 
the following are equivalent.
\begin{enumerate}
\item $Z\subset X$ is topologically normal  for every subscheme of codimension $\geq m$.
\item  For every  point $x\in X$ of codimension $\geq m$ any of the following holds:
\begin{enumerate}
\item  $(x, X)$  is  topologically normal,
\item the   Henselization
 $\bigl(x^{\rm h}, X^{\rm h}\bigr)$  is  topologically normal,
\item   the  strict Henselization
 $\bigl(x^{\rm sh}, X^{\rm sh}\bigr)$  is  topologically normal,
\item    $X^{\rm sh}\setminus\{x^{\rm sh}\}$ is connected or
\item   the  completion
 $\bigl(\hat x, \hat X\bigr)$  is  topologically normal. 
\end{enumerate}
\item  For every quasi-finite \'etale morphism  $g:Y\to X$ with $Y$  connected,
the complement of any closed subset $W\subset Y$ of codimension $\geq m$
 is connected.
\end{enumerate}
\end{thm}

Proof. The equivalence of (1) and (2.a) follows by noting that the
arguments of Proposition \ref{local.nature.prop} preserve the codimension.

(2.a) is equivalent to (2.b) and (2.c) by 
(\ref{hart.ff.up.lem}) and Proposition \ref{reg.morph.preserves.P} and to
(2.e) by (\ref{formal.nature.thm}). (2.c) and (2.d) are equivalent by (\ref{hesnelian.topext.conn.prop}). 

Finally (\ref{topext.conn.prop}) shows that  (2.c) $\Rightarrow$ (3);
here we also use that quasifinite morphisms  do not decrease the codimension.
Since the  strict Henselization is the limit of
 quasi-finite \'etale morphisms, (3) $\Rightarrow$ (2.d) is clear. \qed

\begin{lem}\label{topext.conn.prop}
If $Z\subset X$ is topologically normal then
$X$ and $X\setminus Z$ have the same  connected components.
\end{lem}

Proof.  We may assume that $X$ is connected.
If $X\setminus Z$ is disconnected, write it as $Y_1\cup Y_2$ where the $Y_i$ are disjoint, closed subschemes.
Let $\bar Y_i\subset X$ denote the closure of $Y_i$.
Then the natural map $\bar Y_1\amalg \bar Y_2\to X$ shows that
$Z\subset X$ is not a topologically normal pair. \qed

\begin{prop}\label{hesnelian.topext.conn.prop}
Let $(x, X)$ be a Henselian local scheme with
separably closed residue field $k(x)$.   The following are equivalent.
\begin{enumerate}
\item $(x, X)$ is strongly topologically normal,
\item $(x, X)$  is topologically normal and
\item $X\setminus \{x\}$ is connected.
\end{enumerate}
\end{prop}

Proof.  (1) $\Rightarrow$ (2) is clear and
(2) $\Rightarrow$ (3) follows from (\ref{topext.conn.prop}).
To see (3) $\Rightarrow$ (1), let 
 $g:X'\to X$ be a finite morphism that is a universal homeomorphism
 over $X\setminus\{x\}$. 
If $X$ is Henselian, 
there is a one--to--one correspondence between
connected components of $X'$ and connected components of $x':=g^{-1}(x)$; see
\cite[I.4.2]{milne}. If 
  $k(x')/k(x)$ is not purely inseparable then 
$x'$  is not connected hence
 $X'\setminus\{x'\}$ is also not connected.
The latter is, however, homeomorphic to 
$X\setminus\{x\}$. \qed

\section{Nagata schemes}

\begin{defn} \label{nagata.defn} An integral 
  Noetherian scheme $X$ is called {\it N-1} if its normalization in $k(X)$
is finite over $X$ and  {\it N-2} if its normalization in any finite field extension of $k(X)$
is finite over $X$. The latter is equivalent to the following:
every integral scheme  with a finite, dominant morphism  $X'\to X$ is N-1.

A  Noetherian scheme $X$ is {\it Nagata} if   every integral subscheme
$W\subset X$ is N-2; see for instance  
\cite[p.264]{mats-cr} or \cite[Tag 033R]{stacks-project}.
\end{defn}

We obtain a characterization of local Nagata
schemes in terms of punctual normalizations and  formal fibers. Most parts of this have been known.
A local Nagata domain is analytically unramified and an  analytically unramified domain is N-1; see
\cite[32.2]{MR0155856} and  \cite[Tag 0331]{stacks-project}.
(Recall that  a reduced, semilocal  scheme $(x,X)$ is called
{\it  analytically unramified} iff its completion $\hat X$ is reduced.)
 The connection with formal fibers is
essentially  part of the argument that shows that a
quasi-excellent scheme  is Nagata, see  \cite[33.H]{mats-ca} or \cite[Tag 07QV]{stacks-project}. 
\smallskip

For a  scheme $X$ let $\operatorname{SiS}(X)$ denote the
set of all semilocal, integral  $X$-schemes 
obtained by localizing an  $X$-scheme that is  
 finite over $X$.

\begin{thm} \label{nagata.char.thm}
For a Noetherian scheme the following are equivalent.
\begin{enumerate}
\item All local rings of $X$ are Nagata.
\item For every $(y, Y)\in \operatorname{SiS}(X)$ the normalization  
 $ Y^{\rm n}\to  Y$ is finite.
\item For every $(y, Y)\in \operatorname{SiS}(X)$ the
punctual normalization  $\bigl(y^{\rm pn}, Y^{\rm pn}\bigr)\to (y, Y)$ is finite.
\item Every $(y, Y)\in \operatorname{SiS}(X)$ is analytically unramified.
\item For every $(y, Y)\in \operatorname{SiS}(X)$ the geometric generic fibers
of $(\hat y, \hat Y)\to (y, Y)$ are reduced. 
\item For every $x\in X$ the morphism $(\hat x, \hat X)\to (x, X)$ has (geometrically)  reduced fibers.
\end{enumerate}
\end{thm}

\begin{rem} \label{punc.nagata.char.rem}
This suggests that one should call a  Noetherian  scheme
{\it locally Nagata} if it satsifies the condition 
(\ref{nagata.char.thm}.2). By the theorem, this is equivalent to assuming that
all local rings of $X$ are Nagata.

The condition (\ref{nagata.char.thm}.5) should be compared with the defining property of $G$-rings: For every $x\in X$ the morphism $(\hat x, \hat X)\to (x, X)$ has (geometrically)  regular fibers; see
\cite[p.256]{mats-cr}    or \cite[Tag 07GG]{stacks-project}.
 Hence every G-ring is
locally Nagata.
\end{rem}

\begin{say}[Proof of Theorem \ref{nagata.char.thm}] We may assume that $X$ is integral and affine.
By  induction we may assume that  the equivalence holds over any  integral $X$-scheme such that $Z\to X$ is finite and not dominant.

First note that (1) $\Rightarrow$ (2) by definition and 
(2) $\Rightarrow$ (3) since $ Y^{\rm n}$ dominates
$\bigl(y^{\rm pn}, Y^{\rm pn}\bigr)$.

Assume that (3) holds and fix a point $x\in X$. We aim to show that
$(\hat x, \hat X)\to (x, X)$ has   reduced geometric generic fibers.
Pick any $g\in m_x$. As in the proof of Theorem \ref{punc.nag.lem} (and
using Corollary \ref{cart.pair.nomal.cor}), after replacing $X$ by a suitable  finite  partial normalization,  we may assume that  $X$ is regular at the generic points of $(g=0)$ and the latter has no embedded primes.
By induction, we already know that  $\red (g=0)$ is analytically unramified.
Under these conditions, \cite[Tag 0330]{stacks-project}
says that $X$ is analytically unramified.  Thus (3) $\Rightarrow$ (4).

In order to see the equivalence of (5) and (6) 
fix $\pi:Y\to X$ and set $x:=\pi(y)$. 
Note  that
$\hat Y=Y\times_X\hat X$ where we complete $Y$ at $y$ and $X$ at $x$. If  $\hat X\to  X$ has geometrically reduced fibers then so does $ \hat Y\to  Y$.  Conversely, pick any point $p\in X$ and let
$L\supset k(p)$ be a finite field extension. We can realize $\spec L\to p$
as the generic fiber of a suitable 
$(y, Y)\in \operatorname{SiS}(X)$. Thus if (5) holds for $Y$ then  $\spec L\times_p \hat X$ is
reduced, hence (6) holds.
(It is sufficient to check that the fibers stay reduced after finite field extensions; see \cite[Tag 030V]{stacks-project}.) 

It is clear that (4) implies (5) and if 
  $\hat X\to  X$ has geometrically reduced fibers then so does $ \hat Y\to  Y$. Thus $\hat Y$ is reduced and  $Y$ is
 analytically unramified, proving (6) $\Rightarrow $ (4).

Finally  a strictly increasing chain of partial normalizations
$\cdots \to Y_{i+1}\to Y_i\to \cdots \to Y$ gives a 
a strictly increasing chain of morphisms
$$
\cdots \to \hat Y\times_YY_{i+1}\to \hat Y\times_YY_i\to \cdots \to \hat Y
$$
and these are partial normalizations if $\hat Y$ is reduced. 
This is impossible since a complete semilocal ring is  Nagata  \cite[32.1]{MR0155856}, thus (4) implies (1). \qed
\end{say}

\section{Universal homeomorphisms}

We recall the basic propeties of universal homeomorphisms that we use.

\begin{defn}\label{univ.homeo.defn}
 A morphism of schemes $g:U\to V$ is a 
{\it universal homeomorphism} if  
for every $W\to V$ the
induced morphism $U\times_VW\to W$ is  a homeomorphism.
This notion is called  ``radiciel'' in \cite[I.3.7]{ega71},
``radicial'' in \cite[Tag 01S2]{stacks-project}
and a ``purely inseparable morphism'' by some authors. 
The reason for the latter is the following observation.

Pick a point $v\in V$ and base-change to the algebraic closure of $k(v)$. 
We obtain that  the set-theoretic fiber  $g^{-1}(v)$ is a single 
point $v'$ and 
$k(v')$ is a purely inseparable field extension of   $k(v)$.
We say that $g$ is {\it residue field preserving} if 
$k(v')=k(v)$ for every $v\in V$.

A universal homeomorphism is affine and integral; see
\cite[Tag 04DF]{stacks-project}.
For integral morphisms  the notion of  universal homeomorphism 
is pretty much set theoretic
since a continuous proper map of topological spaces
which is injective and surjective
is a homeomorphism.

\end{defn}

 The following summarizes the basic characterizations.

\begin{lem} \cite[I.3.7--8]{ega71} \label{unihom.char.lem}
 For an integral morphism 
 $g:U\to V$ the following are equivalent.
\begin{enumerate}
\item $g$ is  a universal homeomorphism.
\item $g$ is surjective and universally injective.
\item   $k\bigl(\red g^{-1}(v)\bigr)/k(v)$ is a purely inseparable for every $v\in V$.
\item $g$ is surjective and injective on geometric points.
\end{enumerate}
\end{lem}

In concrete situations it is usually easiest to check (\ref{unihom.char.lem}.3).
Note that for a finite type morphism  $g:U\to V$ the
 points $v\in V$ that satisfy (\ref{unihom.char.lem}.3) form a
constructible set. 

Using property (\ref{unihom.char.lem}.3) we also obtain that if $p:V'\to V$ is
surjective then an integral morphism $g:U\to V$ is a
 universal homeomorphism iff the pull back
$p^*g:V'\times_VU\to V'$ is.

\section{Weakly normal and seminormal schemes}

We recall the definitions and basic properties of 
weakly normal and seminormal schemes, following
\cite{MR0239118, MR0266923, MR571055, MR556308}.

\begin{defn}
\label{weak.norm.defn} 
A morphism of schemes
$g:X'\to X$ is called a {\it partial normalization}
if  $X'$ is reduced, $g$ is integral and $\red(g):X'\to \red(X)$ is 
 birational \cite[1.11]{kk-singbook}. 
Thus a partial  normalization is  dominated by  the normalization $(\red X)^{\rm n}$
of $\red X$ which is the limit of all finite, partial normalizations.
A scheme is normal iff every 
finite, partial normalization $g:X'\to X$ is an isomorphism.

A partial normalization $g:X'\to X$ is
 called a {\it partial seminormalization}
if, in addition, $k\bigl(\red g^{-1}(x)\bigr)=g^*k(x)$
  for every  $x\in X$ 
and a 
{\it partial weak normalization}
if $k\bigl(\red g^{-1}(x)\bigr)/g^*k(x)$ is a purely inseparable field extension
 for every  $x\in X$.
(If $X$ has residue characteristic 0 then 
these 2 notions  coincide.)

There is a
unique largest
partial weak  (resp.\ semi)  normalization  $X^{\rm wn}\to X$ 
(resp.\  $X^{\rm sn}\to X$ ) that dominates
 every other partial weak  (resp.\ semi)   normalization of $X$. It is
 called the 
  {\it  weak normalization}
(resp.\  {\it seminormalization}) of $X$. 

The weak   normalization  $X^{\rm wn}\to X$ is
a universal homeomorphism and it is the maximal
universal homeomorphism that is dominated by the normalization.

A scheme $X$ is called {\it weakly normal} (resp.\  {\it seminormal})
iff every finite, 
partial  weak  (resp.\ semi)  normalization $g:X'\to X$ is an isomorphism.

Note that seminormalization is a functor; that is, every
morphism  $f:Y\to X$ lifts to the seminormalizations
$f^{\rm sn}:Y^{\rm sn}\to X^{\rm sn}$. By contrast,
$f:Y\to X$ lifts to the weak normalizations
$f^{\rm wn}:Y^{\rm wn}\to X^{\rm wn}$ if every irreducible component of $Y$
dominates an irreducible component of $X$ but not otherwise.

 It is easy to see using (\ref{extension.prop}) that an open subscheme
of a weakly  (resp.\ semi) normal scheme is also weakly 
 (resp.\ semi) normal and 
 being  weakly  (resp.\ semi) normal is a local property.
For excellent schemes
  weak  (resp.\ semi) normality   a formal-local property; see
\cite{MR571055, MR556308}.
\end{defn}

\begin{defn}[Constrained normalization]\label{weak.norm.2.defn} 
More generally, let $\pi:Y\to X$ be any  morphism. The
{\it $\pi$-constrained normalization}  
(resp.\ {\it weak  or seminormalization}) of $X$ 
 is the
unique largest partial normalization  (resp.\ weak  or semi  normalization)
$X'\to X$ such that $\red(\pi)$ factors as $\red(Y)\to X'\to X$. 
(Note that we do not require $\pi$ to be dominant.)

The weak  or semi  normalization constrained by the relative normalization
$X_Z^{\rm rn}\to X$ is called the  
{\it relative weak  or seminormalization} of the pair $Z\subset X$.
These are denoted by $X_Z^{\rm rwn}\to X$ and $X_Z^{\rm rsn}\to X$.

As in \cite[Tag 0BAK]{stacks-project}, let $A^{\rm int}\subset \pi_*\o_Y$ 
denote the maximal subalgebra that is integral over $\o_X$ and
 $A^{\rm pi}\subset \pi_*\o_Y$ 
 the maximal subalgebra that is purely inseparable over $\o_X$.
(That is,  $\spec_XA^{\rm pi}\to X$ is a universal homeomorphism and maximal
with this property.)
Then  $\spec_XA^{\rm int}$ is called the {\it relative normalization}
of $X$ in $Y$ and  $\spec_XA^{\rm pi}$
the {\it purely inseparable closure} of $X$ in $Y$.
(The relative normalization frequently agrees with the 
constrained normalization but the notions are quite different if
$\pi$ is not birational.)

\end{defn}

\begin{exmp}\label{weak.closure.field case.exmp}
Let $k$ be a field and $i:k\into A$ a finite, semisimple $k$-algebra;
that is, a sum of finite field extensions. 
If $A$ is separable
then $k$ is purely inseparably closed in $A$ but the converse does not hold.
For instance, $k(s^p, t^p)$ is purely inseparably closed in
$A:=k(s, t^p)+k(s^p, t)$ where $p=\chr k$. 

In order to get a  characterization,
\cite{MR556308} considers the  complex
$$
k\stackrel{(i,i)}{\longrightarrow} A+A\stackrel{\delta}{\longrightarrow}
A\otimes_kA
\eqno{(\ref{weak.closure.field case.exmp}.1)}
$$
 where $\delta(a_1,a_2)=a_1\otimes 1-1\otimes a_2$. 
Let $\red(A\otimes_kA )$ denote the quotient of $A\otimes_kA$ by its
nil-radical. We claim that
$$
k\stackrel{(i,i)}{\longrightarrow} A+A\stackrel{\delta}{\longrightarrow}
\red(A\otimes_kA)
\eqno{(\ref{weak.closure.field case.exmp}.2)}
$$
is exact iff $\spec k$ is weakly normal in $\spec A$.
To see this, pick $(a_1, a_2)\in A+A$ and set $q=(\chr k)^m$ for some $m\geq 1$.
Then 
$$
(a_1\otimes 1-1\otimes a_2)^q=(a_1\otimes 1)^q-(1\otimes a_2)^q=
a_1^q\otimes 1-1\otimes a_2^q, 
$$
hence $\delta(a_1, a_2)$ is nilpotent iff  $a_1^q=a_2^q\in k$ for some
$q=(\chr k)^m$.

Let $Y_k$ be a geometrically regular $k$-scheme.   Then 
$$
\red(A\otimes_kA)\otimes_k\o_{Y_k}=\red\bigl(\o_{Y_A}\otimes_{\o_{Y_k}}\o_{Y_A}\bigr),
$$
thus, if the sequence (\ref{weak.closure.field case.exmp}.2) is exact,
then so is the tensored  sequence 
$$
\o_X\stackrel{(i,i)}{\longrightarrow} \o_{Y_A}+\o_{Y_A}\stackrel{\delta}{\longrightarrow}
\red\bigl(\o_{Y_A}\otimes_{\o_{Y_k}}\o_{Y_A}\bigr).
\eqno{(\ref{weak.closure.field case.exmp}.3)}
$$
Therefore, if $k$ is purely inseparably closed in  $A$ then
 $Y_k$ is purely inseparably closed in $Y_A$.
\end{exmp}

The next result  connecting the absolute and relative versions of
weak and seminormality is very useful in  inductive treatments.
The proof follows \cite{MR571055, MR556308}.

\begin{prop} \label{wn.sn.conductor.crit}
Let $f:Y\to X$ be a finite modification
with conductor subschemes  $C_Y\subset Y$ and
$C_X\subset X$.  
Assume that the pair $C_Y\subset Y$ satisfies WN (resp.\ SN).

Then  the pair $C_X\subset X$ satisfies WN (resp.\ SN)
iff  $C_Y$ (and hence also $C_X$) are  reduced and $ C_X$ is   its own
purely inseparable closure in $C_Y$
 (resp.\ its own  $f|_{C_Y}$-constrained   seminormalization).
\end{prop}

Proof: If $C_Y$ is not reduced then
$\spec_X\bigl(\o_X+\sqrt{\cond_Y}\bigr)\to X$ is a 
universal homeomorphism that is not an isomorphism,
where $\cond_Y\subset \o_Y$ denotes the conductor ideal.
 Let $C_Y\to \tilde C\to C_X$
denote the purely inseparable closure (resp.\  $f|_{C_Y}$-constrained   
seminormalization)  and
$F\subset f_*\o_Y$ the preimage of $\o_{\tilde C}$. 
Then $\spec_X F\to X$ is a universal homeomorphism
(res.\ also preserves residue fields)   that is not an isomorphism.
Thus the conditions are necessary.

Conversely, assume that  $p:X'\to X$ is a universal homeomorphism
(resp.\ also preserves residue fields) that is an isomorphism over
$X\setminus C_X$. Let $\o_{Y'}$ be the composite of
$f_*\o_Y$ and $p_*\o_{X'}$ inside $j_*\o_U$ where $U:=X\setminus C_X$.
Then $Y'\to Y$ is a universal homeomorphism
(resp.\ also preserves residue fields) hence an isomorphism by
assumption. 
  Thus $p_*\o_{X'}\subset f_*\o_Y$ and hence $p_*\o_{X'}/\cond_X\subset \o_{C_Y}$.
The corresponding map 
$\spec_X\bigl(p_*\o_{X'}/\cond_X\bigr)\to C_X$
is a universal homeomorphism (resp.\ partial seminormalization),
hence an isomorphism if $C_Y$ is reduced and $ C_X$ is   its own
purely inseparable closure in $\red C_Y$
 (resp.\  $f|_{C_Y}$-constrained   seminormalization). 
Thus $p:X'\to X$ is an
isomorphism. \qed

\section{Open problems}

Let $R$ be an integral domain with normalization  $R^{\rm n}$. 
Then $R$ is not punctually normal at a prime ideal ${\mathfrak p}\subset R$
iff ${\mathfrak p}$ is an  associated prime of $R^{\rm n}/R$. 
It is not hard to construct 1-dimensional, Noetherian, seminormal domains
with infinitely many  prime ideals, all of which are  associated primes of $R^{\rm n}/R$.

The situation seems more complicated in higher dimensions.
By (\ref{nonrom.finite.inZ.lem}), for every ${\mathfrak q}\neq (0)$, only
fintely many of the primes ${\mathfrak p}\supset {\mathfrak q}$
can be associated primes of $R^{\rm n}/R$. This leads to the following
problems.

\begin{ques} \label{ques.1}
Are there  integral domains of every dimension such that
every height 1 prime  is an associated prime of $R^{\rm n}/R$?
\end{ques}

\begin{ques} Are there seminormal  integral domains of every dimension such that
every height 1 prime  is an associated prime of $R^{\rm n}/R$?
\end{ques}

\begin{ques} What can one say about the  height $\geq 2$  associated primes of $R^{\rm n}/R$?
\end{ques}


\def\cprime{$'$} \def\cprime{$'$} \def\cprime{$'$} \def\cprime{$'$}
  \def\cprime{$'$} \def\cprime{$'$} \def\dbar{\leavevmode\hbox to
  0pt{\hskip.2ex \accent"16\hss}d} \def\cprime{$'$} \def\cprime{$'$}
  \def\polhk#1{\setbox0=\hbox{#1}{\ooalign{\hidewidth
  \lower1.5ex\hbox{`}\hidewidth\crcr\unhbox0}}} \def\cprime{$'$}
  \def\cprime{$'$} \def\cprime{$'$} \def\cprime{$'$}
  \def\polhk#1{\setbox0=\hbox{#1}{\ooalign{\hidewidth
  \lower1.5ex\hbox{`}\hidewidth\crcr\unhbox0}}} \def\cdprime{$''$}
  \def\cprime{$'$} \def\cprime{$'$} \def\cprime{$'$} \def\cprime{$'$}
\providecommand{\bysame}{\leavevmode\hbox to3em{\hrulefill}\thinspace}
\providecommand{\MR}{\relax\ifhmode\unskip\space\fi MR }
\providecommand{\MRhref}[2]{%
  \href{http://www.ams.org/mathscinet-getitem?mr=#1}{#2}
}
\providecommand{\href}[2]{#2}

\medskip

\noindent JK: Princeton University, Princeton NJ 08544-1000

{\begin{verbatim} kollar@math.princeton.edu\end{verbatim}}

\end{document}